\documentclass[a4paper,leqno]{article}

\usepackage{amsmath}
\usepackage{amssymb}
\usepackage{amsfonts}
\usepackage{enumerate}
\usepackage{theorem}

\theoremstyle{plain}
\newtheorem{teo}{Theorem}[section]
\newtheorem{lem}[teo]{Lemma}

\newtheorem{prop}[teo]{Proposition}

{\theorembodyfont{\rmfamily}

\newtheorem{exam}[teo]{Example}
}

\renewcommand{\eqref}[1]{\textnormal{(\ref{#1})}}

\numberwithin{equation}{section}

\newcommand{\cvd}{\hfill$\square$}
\newcommand{\proof}[1]{\noindent\textsc{Proof#1}}

\title{Continuity properties of Neumann-to-Dirichlet maps with respect to the $H$-convergence of the coefficient matrices}

\author{Luca \textsc{Rondi}\\
\normalsize{Dipartimento di Matematica e Geoscienze}\\
\normalsize{Universit\`a degli Studi di Trieste}\\
\normalsize{via Valerio, 12/1}\\
\normalsize{34127 Trieste, Italy}\\
\normalsize{\texttt{rondi@units.it}}}

\date{}

\begin{document}

\maketitle

\setcounter{section}{0}
\setcounter{secnumdepth}{2}

\begin{abstract}
We investigate the continuity of boundary operators, such as the Neu\-mann-to-Dirichlet map, with respect to the coefficient matrices of the underlying elliptic equations. We show that for nonsmooth coefficients the correct notion of convergence is the one provided by $H$-convergence (or $G$-convergence for symmetric matrices). We prove existence results for minimum problems associated to variational methods used to solve the so-called inverse conductivity problem, at least if we allow the conductivities to be anisotropic. In the case of isotropic conductivities we show that on certain occasions existence of a minimizer may fail.

\medskip


\medskip

\end{abstract}

\section{Introduction}

We assume that a conducting body is contained in a bounded domain $\Omega\subset\mathbb{R}^N$, $N\geq 2$, with Lipschitz boundary $\partial\Omega$. Let the conductivity tensor 
of the body be given by a matrix $A=A(x)$, $x\in\Omega$, satisfying suitable ellipticity conditions.
Given a current density $g$ on the boundary, for instance $g\in L^2(\partial\Omega)$ with zero mean, the electrostatic potential $v$ in $\Omega$ is the unique solution to
\begin{equation}\label{directproblem}
\left\{\begin{array}{ll}
\mathrm{div}(A\nabla v)=0&\text{in }\Omega\\
A\nabla v\cdot\nu=g&\text{on }\partial\Omega\\
\int_{\partial\Omega}v=0
\end{array}\right.
\end{equation}
where $\nu$ denotes the exterior unit normal.
We define the \emph{Neumann-to-Dirichlet map} associated to $A$ the linear and bounded operator $N\text{-}D(A):L^2_{\ast}(\partial\Omega)\to L^2_{\ast}(\partial\Omega)$ such that, for any $g\in L^2_{\ast}(\partial\Omega)$,
$N\text{-}D(A)[g]=v|_{\partial\Omega}$, $v$ solution to \eqref{directproblem}.
Here $L^2_{\ast}(\partial\Omega)$ denotes the space of $L^2(\partial\Omega)$ functions
with zero mean.

The inverse conductivity problem, proposed by Calder\'on in \cite{Cal}, consists of determining an unknown conductivity tensor $A$ from electrostatic measurements at the boundary of voltage and current type, namely from the knowledge of its corresponding Neumann-to-Dirichlet map $N\text{-}D(A)$. Uniqueness has been extensively investigated for the case of scalar, that is isotropic, conductivities. In dimension $3$ and higher, the breakthrough was in the middle 80's, first by 
Kohn and Vogelius for the determination of the conductivity at the boundary and for the
 analytic case, \cite{Koh e Vog84:1,Koh e Vog85}, then by Sylvester and Uhlmann for smooth $C^2$ conductivities, \cite{Syl e Uhl87}. Almost ten years later the two dimensional case for smooth conductivities was solved by Nachman, \cite{Nac}.

Much more recently these uniqueness results have been greatly improved. For $N\geq 3$  the a priori regularity conditions on the unknown conductivity have been relaxed up to Lipschitz, \cite{Hab e Tat}, and even to conductivities whose gradient is unbounded, \cite{Hab}.
For $N=2$ we have uniqueness without any regularity assumptions, \cite{Ast e Pai}. Therefore for $N=2$ the uniqueness issue is completely solved for scalar conductivities.

Since the Neumann-to-Dirichlet map is invariant under suitable changes of coordinates that keep fixed the boundary, uniqueness is never achieved for anisotropic conductivities. However, at least in dimension $2$ for symmetric conductivity tensors, this is the only obstruction as shown first in \cite{Syl} in the smooth case and then in \cite{Ast e Pai e Las} in the general case.

We consider numerical methods of reconstruction, especially those of a variational character.
Let us fix the ideas and consider the following least squares variational approach
\begin{equation}\label{minimizationexample}
\min_{A\in\mathcal{M}}\| N\text{-}D(A)-\hat{\Lambda}\|
\end{equation}
where $\hat{\Lambda}$ is the measured, therefore approximated, Neumann-to-Dirichlet map, $\mathcal{M}$ is a suitable class of admissible conductivity tensors and $\|\cdot\|$ is a suitable operator norm.

In order to obtain existence of a minimizer
it is crucial to establish continuity, or at least lower semicontinuity, properties of the forward operator $N\text{-}D$ that to each conductivity tensor $A$ associates the corresponding Neumann-to-Dirichlet map $N\text{-}D(A)$.
These continuity properties clearly strongly depends on the topology of the space of admissible conductivity tensors and, up to a certain extent, also on the norm used for the Neumann-to-Dirichlet maps. For instance it is an easy remark that if we use the $L^{\infty}$ norm on the conductivity tensors, then $N\text{-}D$ is continuous with respect to the natural norm on the Neumann-to-Dirichlet maps, that is the one of a linear operator between suitable Sobolev spaces on $\partial\Omega$, see Theorem~\ref{inftycontinuity}. This result is of interest only if one deals with smooth conductivities tensors. In fact, for discontinuous conductivities the $L^{\infty}$ norm is not suited, see for instance the inclusion problem. In \cite{Ron08} the continuity with respect to $L^q$ norms, with $q$ finite, was considered and continuity was proved if we consider a slightly weaker topology on the Neumann-to-Dirichlet maps, namely the one of a linear operator between $L^2_{\ast}(\partial\Omega)$ and itself. We recall this result in Theorem~\ref{pcontinuity}. In this introduction we shall refer to this norm as the $L^2\text{-}L^2$ norm, in contrast with the natural norm, on Neumann-to-Dirichlet maps.

However, the main disadvantage for using (strong) $L^q$ convergence, with $q$ finite or not,
is that we do not have compactness, unless we use a regularization method or impose suitable a priori conditions. Let us mention here that in this paper we do not address the issue of ill-posedness and hence of the regularization that is required for the numerical solution of these minimization problems. Due to the lack of compactness, existence of solutions of minimum problems such as \eqref{minimizationexample} may not be guaranteed by using strong $L^q$ convergence. On the other hand, weak convergence in $L^q$ spaces with $q$ finite or weak$\ast$ convergence in $L^{\infty}$ does not lead to convergence of solutions of the corresponding elliptic equations, as it is well-known in homogenization theory. Again inspired by homogenization theory, the correct notion of convergence for conductivity tensors that leads to convergence of solutions of the corresponding elliptic equations is $H$-convergence. The great advantage of $H$-convergence is that it does not require any regularity and it has compactness properties. Notice that $H$-convergence reduces to $G$-convergence for symmetric conductivity tensors. For simplicity in this introduction we limit ourselves to this case and then we consider $G$-convergence only.

We remark that the relationship between $G$-convergence and the Dirichlet-to-Neumann maps has been already studied in \cite{Ale e Cab} and \cite{Far e Kur e Rui}. In \cite{Far e Kur e Rui} a
continuity result is established for Dirichlet-to-Neumann maps, endowed with their natural norm, if we have a sequence of symmetric conductivity tensors $G$-converging and converging in a suitable sense near the boundary to another symmetric conductivity tensor. Then their results are applied to cloaking.

Our interest is in the variational approach to reconstruction. In Section~\ref{applicationsection}, the main of the paper, we analyze several of them that have been proposed in the literature.
First we analyze the case in which all, or at least infinitely many, measurements are used. Then we consider the more practical case of a finite number of measurements,
following the approaches by Yorkey, \cite{Yor}, and Kohn and Vogelius, \cite{Koh e Vog87}, and Kohn and McKenney, \cite{Koh e McK}.
Finally we treat the more realistic case of the so-called experimental measurements.

A key point for proving existence of minimizers is establishing
continuity properties with respect to $G$-convergence of symmetric conductivity tensors. We
obtain a lower semicontinuity result for the Neumann-to-Dirichlet map with respect to the natural norm, see Theorem~\ref{semicontinuitycor}. However, this may be improved to continuity
if we use on the Neumann-to-Dirichlet maps the $L^2\text{-}L^2$ norm, see Theorem~\ref{Hcontinuity}. Notice that we actually may allow nonsymmetric conductivity tensors, that is our result is not limited to symmetric conductivity tensors, and that we do not impose any convergence of the conductivities at the boundary.
Through compactness these results allow to prove existence of a minimizer to \eqref{minimizationexample}, both with the natural norm or with the $L^2\text{-}L^2$ norm,
provided $\mathcal{M}$ is $G$-closed, which is true for classes of symmetric conductivity tensors satisfying fixed ellipticity conditions, see Theorem~\ref{existence}

On the other hand, the issue of solving these minimization problems over a class of scalar conductivities still remains open, because scalar conductivities may actually $G$-converge to an anisotropic symmetric conductivity tensor. Therefore if we use the direct method to solve \eqref{minimizationexample} in a class of scalar conductivities we may not guarantee that a minimizer exists since the $G$-limit of a minimizing sequence might not be a scalar conductivity any more. Of course if this limit is the push-forward of a scalar conductivity by a change of coordinates that keeps fixed the boundary, then this scalar conductivity would be a minimizer. However this may be impossible since an example kindly communicated by Giovanni Alessandrini, Example~\ref{Giovanni}, shows that there are anisotropic symmetric conductivity tensors whose Neumann-to-Dirichlet map is not the Neumann-to-Dirichlet map of any scalar conductivity.

Actually, existence may indeed fail. This is  a nice, or actually not that nice, side effect of our continuity results.
In fact, if we use the $L^2\text{-}L^2$ norm on Neumann-to-Dirichlet maps we may show, Example~\ref{nonexistence}, that for some $\hat{\Lambda}$ a minimizer of \eqref{minimizationexample} in a class of scalar conductivities does not exist. For this nonexistence result, beside continuity,
we make use of  Example~\ref{Giovanni} and of the density of scalar conductivities inside the symmetric conductivity tensors
with respect to $G$-convergence, which is recalled in Proposition~\ref{approxbyscalar}.

We notice that, as pointed out already in Theorem~4.9 in \cite{Far e Kur e Rui}, if we use the Dirichlet-to-Neumann maps with their natural norm we can not expect continuity. However we show that lower semicontinuity holds and this is enough for the variational problems we wish to solve. On the other hand, in the applications it is considered more convenient to prescribe the current density and measure the corresponding potential, that is to use the Neumann-to-Dirichlet map instead of the Dirichlet-to-Neumann map. Moreover, since only a finite number of measurements may actually be performed and the corresponding prescribed current densities may be chosen, it is not restrictive to assume that they have a (mild) regularity namely that they belong to $L^2$. Also the error on the measured voltage may be very naturally defined through an integral norm such as the $L^2$ norm. Therefore, the use of the $L^2\text{-}L^2$ norm for the Neumann-to-Dirichlet maps seems to be the correct choice from the point of view of the applications. Indeed this is actually the variational approach by Yorkey, \cite{Yor}, who, for a finite number of prescribed current densities $g_i$, $i=1,\ldots,n$, considers the following minimization problem
$$
\min_{A\in\mathcal{M}}\sum_{i=1}^n\int_{\partial\Omega}|N\text{-}D(A)[g_i]-\varphi_i|^2
$$
where $\varphi_i$, $i=1,\ldots,n$, are the corresponding measured voltages.
This problem is exactly the finite measurements version of \eqref{minimizationexample} with the $L^2\text{-}L^2$ norm.

Another confirmation on the fact that the $L^2\text{-}L^2$ norm for the Neumann-to-Dirichlet maps is the  right choice comes from the more realistic model of current and voltage measurements, the so-called experimental measurements introduced in \cite{Som e Che e Isa} where the measurements are modeled through a resistance matrix $R$. We show in Proposition~\ref{expmeasprop} that the distance between the resistance matrices corresponding to two different conductivity tensors may be controlled by the $L^2\text{-}L^2$ norm of the difference of the corresponding Neumann-to-Dirichlet maps.

This means that in the applications using the $L^2\text{-}L^2$ norm for the Neumann-to-Dirichlet maps is not a restriction and we obtain good continuity properties on the space of symmetric conductivity tensors with respect to $L^q$ convergence, with $q$ possibly finite, or with respect to $G$-convergence. The same continuity properties clearly hold for the experimental measurements case.

Let us mention here that, in order to treat also the approach proposed 
by Kohn and Vogelius, \cite{Koh e Vog87}, and Kohn and McKenney, \cite{Koh e McK},
we analyze the continuity properties with respect to $G$-convergence of suitable operators
$KN_1$ and $KN_2$, introduced at the end of Section~\ref{notationsection}. Similar results and considerations hold also for these operators, in particular for $KN_1$ the more interesting of the two.

%

The plan of the paper is the following. In the preliminary Section~\ref{notationsection} we set the notation, define the operators of our interest, and review the definition and basic properties of $H$-convergence. In Section~\ref{applicationsection} we discuss existence and nonexistence issues for variational methods of reconstruction for the inverse conductivity problems proposed in the literature. Section~\ref{Hsection}, the technical core of the paper, contains the proofs of the continuity properties with respect to $H$- and $G$-convergence.
Finally, in the Appendix~\ref{Lpsection} we briefly recall, for the sake of completeness, the continuity properties with respect to $L^q$ norms on the space of conductivity tensors.

\section{Preliminaries}\label{notationsection}

Throughout the paper we shall keep fixed 
positive constants $\alpha$, $\beta$ and $\tilde{\beta}$, with $0<\alpha\leq  \beta,\tilde{\beta}$. For $N\geq 2$, we call $\mathbb{M}^{N\times N}(\mathbb{R})$ the space of real valued $N\times N$ matrices.

Let us consider the following ellipticity conditions for a given $A\in \mathbb{M}^{N\times N}(\mathbb{R})$, with $N\geq 2$. The first condition is the usual one
\begin{equation}\label{ell1}
\left\{\begin{array}{ll}
A\xi\cdot\xi\geq \alpha\|\xi\|^2&\text{ for any }\xi\in\mathbb{R}^N\\
\|A\|\leq \beta
\end{array}\right.
\end{equation}
where, for any $N\times N$ matrix $A$, $\|A\|$ denotes its norm as a linear operator of $\mathbb{R}^N$ into itself.
The second ellipticity condition is the following
\begin{equation}\label{ell2}
\left\{\begin{array}{ll}
A\xi\cdot\xi\geq \alpha\|\xi\|^2&\text{ for any }\xi\in\mathbb{R}^N\\
A^{-1}\xi\cdot\xi\geq \tilde{\beta}^{-1}\|\xi\|^2&\text{ for any }\xi\in\mathbb{R}^N.\end{array}\right.
\end{equation}
Let us notice that if $A$ satisfies \eqref{ell2} with constants $\alpha$ and $\tilde{\beta}$, then it also satisfies \eqref{ell1} with constants $\alpha$ and $\beta=\tilde{\beta}$. On the other hand, if $A$ satisfies \eqref{ell1} with constants $\alpha$ and $\beta$, then it also satisfies \eqref{ell2} with constants $\alpha$ and $\tilde{\beta}=\beta^2/\alpha$. If $A$ is symmetric then \eqref{ell1} and \eqref{ell2} are equivalent and both correspond to the condition, with $\beta=\tilde{\beta}$,
$$\alpha\|\xi\|^2\leq A\xi\cdot\xi\leq \beta\|\xi\|^2\quad \text{ for any }\xi\in\mathbb{R}^N.$$
Finally, if $A=\sigma I_N$, where $I_N$ is the $N\times N$ identity matrix and $\sigma$ is a real number, the condition further reduces to
$$\alpha\leq \sigma\leq \beta.$$

We define the following classes of conductivity tensors in $\Omega$,
$\Omega\subset\mathbb{R}^N$ being a bounded open set. We call
$\mathcal{M}(\alpha,\beta)$, respectively $\tilde{\mathcal{M}}(\alpha,\tilde{\beta})$, the set of $A=A(x)$, $x\in \Omega$, an $N\times N$ matrix whose entries are real valued  measurable functions in $\Omega$, such that, for almost any $x\in\Omega$, $A(x)$ satisfies \eqref{ell1}, respectively \eqref{ell2}. We call $\mathcal{M}_{sym}(\alpha,\beta)$, respectively $\mathcal{M}_{scal}(\alpha,\beta)$, the set of $A\in \mathcal{M}(\alpha,\beta)$ such that, for almost any $x\in\Omega$, $A(x)$ is symmetric, respectively  $A(x)=\sigma(x) I_N$ with $\sigma(x)$ a real number.
Obviously we have $\tilde{\mathcal{M}}(\alpha,\tilde{\beta})\subset \mathcal{M}(\alpha,\tilde{\beta})$
and $\mathcal{M}(\alpha,\beta)\subset \tilde{\mathcal{M}}(\alpha,\beta^2/\alpha)$.
We recall that by a conductivity tensor $A$ in $\Omega$, respectively symmetric conductivity tensor or scalar conductivity, we mean $A\in \mathcal{M}(\alpha,\beta)$, respectively
$\mathcal{M}_{sym}(\alpha,\beta)$ or $\mathcal{M}_{scal}(\alpha,\beta)$, for some constants $0<\alpha\leq\beta$.

We notice that all these classes are closed with respect to the $L^p$ metric, for any $p$, $1\leq p\leq +\infty$, where for any two conductivity tensors $A_1$ and $A_2$ in $\Omega$
$$\|A_1-A_2\|_{L^p(\Omega)}=\|(\|A_1-A_2\|)\|_{L^p(\Omega)}.$$

We wish to investigate the continuity of solutions to Dirichlet and Neumann boundary value problems with respect to the conductivity tensor $A$. We begin with the Dirichlet case.

Let $\Omega\subset \mathbb{R}^N$, $N\geq 2$, be a bounded open set.
Let us fix $\varphi\in W^{1,2}(\Omega)$ and $F\in (W^{1,2}(\Omega))^{\ast}$. For any conductivity tensor $A$, let $u=u(A,F,\varphi)$ be the unique weak solution to the following Dirichlet boundary value problem
\begin{equation}\label{Dirichletproblem}
\left\{\begin{array}{ll}
-\mathrm{div}(A\nabla u)=F&\text{in }\Omega\\
u=\varphi&\text{on }\partial\Omega.
\end{array}\right. 
\end{equation}

The weak formulation is to look for $u\in W^{1,2}(\Omega)$ such that $u-\varphi\in
W^{1,2}_0(\Omega)$ and
$$\int_{\Omega}A\nabla u\cdot\nabla\psi=\langle F,\psi\rangle_{((W^{1,2}(\Omega))^{\ast},W^{1,2}(\Omega))}\quad\text{for any }\psi\in W_0^{1,2}(\Omega).$$


Let us consider two conductivity tensors $A_1$ and $A_2\in\tilde{\mathcal{M}}(\alpha,\tilde{\beta})$.
We denote $u_1=u(A_1,F,\varphi)$ and $u_2=u(A_2,F,\varphi)$.

First of all it is easy to show that
$$\|u_1-u_2\|_{W^{1,2}(\Omega)}\leq C(\|\varphi\|_{W^{1,2}(\Omega)}+\|F\|_{(W^{1,2}\Omega))^{\ast}})\|A_1-A_2\|_{L^{\infty}(\Omega)}$$
where $C$ depends on $N$, $\Omega$, $\alpha$ and $\tilde{\beta}$ only.

Similar reasonings hold if we consider Neumann problems instead of Dirichlet ones.
In this case we need to restrict ourselves to the case in which $\Omega\subset \mathbb{R}^N$, $N\geq 2$, is a bounded domain with Lipschitz boundary.

Let us fix $F\in (W^{1,2}(\Omega))^{\ast}$ such that $\langle F,1\rangle_{((W^{1,2}(\Omega))^{\ast},W^{1,2}(\Omega))}=0$.
For any conductivity tensor $A$, let $v=v(A,F)$ be the unique weak solution to the following Neumann boundary value problem
\begin{equation}\label{Neumannproblem}
\left\{\begin{array}{ll}
-\mathrm{div}(A\nabla v)=F&\text{in }\Omega\\
\int_{\partial\Omega}v=0.
\end{array}\right. 
\end{equation}

The weak formulation is to look for $v\in W^{1,2}(\Omega)$ such that $\int_{\partial\Omega}v=0$ and
\begin{equation}
\int_{\Omega}A\nabla v\cdot\nabla\psi=\langle F,\psi\rangle_{((W^{1,2}(\Omega))^{\ast},W^{1,2}(\Omega))}\quad\text{for any }\psi\in W^{1,2}(\Omega).
\end{equation}

We notice that we use a normalization condition which is suited for our application, namely the one of zero mean on the boundary.
However all the result we shall state for Neumann problems  still hold if we replace this normalization condition with the following more common one
$$\int_{\Omega}v=0.$$

Let us again consider two conductivity tensors $A_1$ and $A_2\in \tilde{\mathcal{M}}(\alpha,\tilde{\beta})$. We denote $v_1=v(A_1,F)$ and $v_2=v(A_2,F)$.

Again it is easy to show that
$$\|v_1-v_2\|_{W^{1,2}(\Omega)}\leq C\|F\|_{(W^{1,2}\Omega))^{\ast}}\|A_1-A_2\|_{L^{\infty}(\Omega)}$$
where $C$ depends on $N$, $\Omega$, $\alpha$ and $\tilde{\beta}$ only.

\subsubsection{Definition of our forward operators}

Let $\Omega\subset \mathbb{R}^N$, $N\geq 2$, be a bounded domain with Lipschitz boundary.
We recall that $W^{1/2,2}(\partial \Omega)$ is the space of traces of
$W^{1,2}(\Omega)$ functions on $\partial \Omega$ and that $W^{-1/2,2}(\partial\Omega)$ is its dual.
Let us notice that $W^{-1/2,2}(\partial\Omega)\subset (W^{1,2}(\Omega))^{\ast}$
by identifying any $g\in W^{-1/2,2}(\partial\Omega)$ with $F(g)$ such that for any $\psi\in W^{1,2}(\Omega)$
$$\langle F(g),\psi\rangle_{((W^{1,2}(\partial \Omega))^{\ast},W^{1,2}(\partial \Omega))}= \ \langle g,\psi|_{\partial\Omega}\rangle_{(W^{-1/2,2}(\partial \Omega),W^{1/2,2}(\partial \Omega))}.$$
Moreover, we call  
$W^{-1/2,2}_{\ast}(\partial\Omega)$ the subspace of 
$g\in W^{-1/2,2}(\partial\Omega)$ such that
$\langle F(g),1\rangle=0$.
We call $W^{1/2,2}_{\ast}(\partial\Omega)$ and $L^2_{\ast}(\partial\Omega)$ the subspaces of $W^{1/2,2}(\partial\Omega)$ and $L^2(\partial\Omega)$ functions with zero mean, respectively. 
Let us notice that $W^{1/2,2}_{\ast}(\partial\Omega)\subset L^2_{\ast}(\partial\Omega)$ and 
$L^2_{\ast}(\partial\Omega)\subset W^{-1/2,2}_{\ast}(\partial\Omega)$ with compact immersions, and $L^2_{\ast}(\partial\Omega)$ is dense in $W^{-1/2,2}_{\ast}(\partial\Omega)$.

We recall that for any two Banach spaces $B_1$, $B_2$, $\mathcal{L}(B_1,B_2)$
will denote the Banach space of bounded linear operators from $B_1$ to $B_2$ with the usual operator norm.

Fixed a conductivity tensor $A$ in $\Omega$, let us define its corresponding Dirichlet-to-Neumann map
$$D\text{-}N(A):W^{1/2,2}(\partial\Omega)\to W^{-1/2,2}(\partial\Omega)$$
where for each $\varphi\in W^{1/2,2}(\partial\Omega)$,
$$D\text{-}N(A)(\varphi)[\psi]=\int_{\Omega}A\nabla u\cdot\nabla \tilde{\psi} \quad\text{for any }\psi\in W^{1/2,2}(\partial\Omega)$$
with $u=u(A,0,\varphi)$ solution to
\begin{equation}\label{Dirichletproblem1}
\left\{\begin{array}{ll}
-\mathrm{div}(A\nabla u)=0&\text{in }\Omega\\
u=\varphi&\text{on }\partial\Omega
\end{array}\right. 
\end{equation}
and $\tilde{\psi}\in W^{1,2}(\Omega)$ is such that
$\tilde{\psi}=\psi$ on $\partial\Omega$ in the trace sense.
We have that $D\text{-}N(A)$ is a well-defined bounded linear operator. Moreover, provided $A\in\tilde{\mathcal{M}}(\alpha,\tilde{\beta})$, its norm is bounded by a constant depending on $N$, $\Omega$, $\alpha$ and $\tilde{\beta}$ only.
Let us notice that, actually, we have $D\text{-}N(A):W^{1/2,2}(\partial\Omega)\to W^{-1/2,2}_{\ast}(\partial\Omega)$.

In an analogous way we define the Neumann-to-Dirichlet map. 
Fixed a conductivity tensor $A$ in $\Omega$, we define its corresponding Neumann-to-Dirichlet map
$$N\text{-}D(A):W^{-1/2,2}_{\ast}(\partial\Omega)\to W^{1/2,2}_{\ast}(\partial\Omega)$$
where for each $g\in W^{-1/2,2}_{\ast}(\partial\Omega)$,
$$N\text{-}D(A)(g)=v|_{\partial\Omega}$$
with $v$ solution to
\begin{equation}\label{Neumannproblem1}
\left\{\begin{array}{ll}
-\mathrm{div}(A\nabla v)=0&\text{in }\Omega\\
A\nabla v\cdot \nu=g&\text{on }\partial\Omega\\
\int_{\partial \Omega}v=0.&
\end{array}\right. 
\end{equation}
Its weak formulation is the following. We look for 
$v\in W^{1,2}(\Omega)$ such that $\int_{\partial\Omega}v=0$ and
$$
\int_{\Omega}A\nabla v\cdot\nabla\psi=\langle g,\psi|_{\partial\Omega}\rangle_{((W^{-1/2,2}(\partial\Omega))^{\ast},W^{1/2,2}(\partial\Omega))}\quad\text{for any }\psi\in W^{1,2}(\Omega).
$$
That is $v=v(A,F(g))$ with $F(g)$ as above. For simplicity, by identifying $F(g)$ with $g$ we shall simply write $v=v(A,g)$.

We have that $N\text{-}D(A)$ is a well-defined bounded linear operator, which is essentially the inverse of $D\text{-}N(A)|_{W^{1/2,2}_{\ast}(\partial\Omega)}$. Moreover, provided $A\in\tilde{\mathcal{M}}(\alpha,\tilde{\beta})$, its norm is bounded by a constant depending on $N$, $\Omega$ and $\alpha$ only.

Our forward operators are 
$$
D\text{-}N:\tilde{\mathcal{M}}(\alpha,\tilde{\beta})
\to \mathcal{L}(W^{1/2,2}(\partial \Omega),W^{-1/2,2}_{\ast}(\partial \Omega))$$ or
$$N\text{-}D:\tilde{\mathcal{M}}(\alpha,\tilde{\beta})
\to \mathcal{L}(W^{-1/2,2}_{\ast}(\partial \Omega),W^{1/2,2}_{\ast}(\partial \Omega)).$$

We shall also be interested in the following operators.

Let us fix 
$\hat{\Lambda}_D\in  \mathcal{L}(W^{1/2,2}(\partial \Omega),W^{-1/2,2}_{\ast}(\partial \Omega))$.
For any conductivity tensor $A$ in $\Omega$, let us define the linear operators
$KD_i(A):W^{1/2,2}(\partial \Omega)\to L^2(\Omega,\mathbb{R}^N)$, $i=1,2$
in the following way. For any $\varphi\in W^{1/2,2}(\partial \Omega)$, we set
$$KD_i(A)[\varphi]=A^{i/2}(\nabla u(A,0,\varphi)-\nabla v(A,\hat{\Lambda}_D(\varphi))).$$

Analogously, if $\hat{\Lambda}_N\in  \mathcal{L}(W^{-1/2,2}_{\ast}(\partial \Omega),W^{1/2,2}_{\ast}(\partial \Omega))$, for any conductivity tensor $A$ in $\Omega$, we define the linear operators
$KN_i(A):W^{-1/2,2}_{\ast}(\partial \Omega)\to L^2(\Omega,\mathbb{R}^N)$, $i=1,2$
in the following way. For any $g\in W^{-1/2,2}_{\ast}(\partial \Omega)$, we set
$$KN_i(A)[g]=A^{i/2}(\nabla u(A,0,\hat{\Lambda}_N(g))-\nabla v(A,g)).$$

We observe that all these are well-defined bounded linear operators. Notice also that, in both cases, for $i=1$ we limit ourselves to symmetric conductivity tensors.

Therefore we also have forward operators
$$
\begin{array}{l}
KD_1:\mathcal{M}_{sym}(\alpha,\beta)
\to \mathcal{L}(W^{1/2,2}(\partial \Omega),L^2(\Omega,\mathbb{R}^N)),\\
KD_2: \tilde{\mathcal{M}}(\alpha,\tilde{\beta})
\to \mathcal{L}(W^{1/2,2}(\partial \Omega),L^2(\Omega,\mathbb{R}^N)),
\end{array}
$$
or
$$
\begin{array}{l}
KN_1:\mathcal{M}_{sym}(\alpha,\beta)
\to \mathcal{L}(W^{-1/2,2}_{\ast}(\partial \Omega),L^2(\Omega,\mathbb{R}^N)),\\
KN_2: \tilde{\mathcal{M}}(\alpha,\tilde{\beta})
\to \mathcal{L}(W^{-1/2,2}_{\ast}(\partial \Omega),L^2(\Omega,\mathbb{R}^N)).
\end{array}
$$

\subsubsection{$H$-convergence}

We recall the definition of $H$-convergence. For its basic properties we refer for instance to \cite{Mur e Tar1},
 \cite{All} and \cite{Cio e Don}.

Let $\Omega\subset\mathbb{R}^N$, $N\geq 2$, be a bounded open set.
For any conductivity tensor $A$ in $\Omega$ and any $F\in (W^{1,2}(\Omega))^{\ast}$, we call $u=u(A,F,0)$ the unique weak solution to the following Dirichlet boundary value problem
$$\left\{\begin{array}{ll}
-\mathrm{div}(A\nabla u)=F&\text{in }\Omega\\
u=0&\text{on }\partial\Omega.
\end{array}\right. 
$$
Given a sequence of conductivity tensors
$A_n\in \tilde{\mathcal{M}}(\alpha,\tilde{\beta})$, $n\in\mathbb{N}$,  we say that $A_n$ $H$-converges to a conductivity tensor $A$, as $n\to\infty$, if for any $F\in (W^{1,2}(\Omega))^{\ast}$ we have
\begin{equation}\label{Hconv}
\left\{\begin{array}{l}
u_n\rightharpoonup u\quad \text{weakly in }W^{1,2}_0(\Omega)\\
A_n\nabla u_n \rightharpoonup A\nabla u \quad \text{weakly in }L^2(\Omega,\mathbb{R}^N)
\end{array}\right. 
\end{equation}
where $u_n=u(A_n,F,0)$, for any $n\in\mathbb{N}$, and $u=u(A,F,0)$.

We make the following observations. First, since $L^2(\Omega)$ is a dense subset of 
$(W^{1,2}(\Omega))^{\ast}$, it is enough that \eqref{Hconv} holds for any $F\in L^2(\Omega)$.
Another important remark is that $\tilde{\mathcal{M}}(\alpha,\tilde{\beta})$ is closed with respect to $H$-convergence, that is the $H$-limit $A$ belongs to $\tilde{\mathcal{M}}(\alpha,\tilde{\beta})$ as well. Finally, if $A_n\in \mathcal{M}_{sym}(\alpha,\beta)$, for any $n\in\mathbb{N}$, and $A_n$ $H$-converges to $A$ as $n\to \infty$, then $A\in \mathcal{M}_{sym}(\alpha,\beta)$ as well. In this case $H$-convergence coincides with the well-known $G$-convergence.

The crucial property of $H$-convergence is the following compactness result.

\begin{teo}
The class $\tilde{\mathcal{M}}(\alpha,\tilde{\beta})$ is \textnormal{(}sequentially\textnormal{)} compact with respect to $H$-convergence. As a corollary, $\mathcal{M}_{sym}(\alpha,\beta)$ is \textnormal{(}sequentially\textnormal{)} compact with respect to $H$-convergence or $G$-convergence.
\end{teo}

Notice, however, that neither $\mathcal{M}(\alpha,\beta)$ nor, more importantly, $\mathcal{M}_{scal}(\alpha,\beta)$ are closed under $H$-convergence. Actually the following proposition
shows that any symmetric conductivity tensor is the $H$-limit of a suitable sequence of scalar conductivities. This is a simple consequence of the characterization of $H$-limits, or better $G$-limits, of sequences of scalar conductivities assuming only two given values, see for instance
Proposition~10 in \cite{Mur e Tar2}, investigated with the purpose of studying the homogenization of composite materials.

\begin{prop}\label{approxbyscalar}
Let $A\in \mathcal{M}_{sym}(\alpha,\beta)$. Then there exists a constant $\alpha_1$, $0<\alpha_1<1$. depending on $\alpha$, $\beta$ and $N$ only, such that, if we call $\beta_1=1/\alpha_1$, the following property holds.

There exists a sequence $A_n\in \mathcal{M}_{scal}(\alpha_1,\beta_1)$, $n\in\mathbb{N}$, such that $A_n$ $H$-converges to $A$ as $n\to \infty$ and for any $n\in\mathbb{N}$ we have that
$A_n=\sigma_n I_N$ where $\sigma_n$ takes only the two values $\alpha_1$ and $\beta_1$.
\end{prop}

\proof{.} 
We fix $\varepsilon$, $0<\varepsilon<1$, and $\theta\in (0,1)$. We consider
$\alpha_1=\varepsilon$ and $\beta_1=1/\varepsilon$. To approximate $A$ by scalar conductivities assuming only values $\alpha_1$ and $\beta_1$ it is enough to show, by \cite[Proposition~10]{Mur e Tar2}, that, for almost every $x\in \Omega$, we have
\begin{equation}\label{suffapprox}
\left\{
\begin{array}{l}
\lambda_-\leq \alpha\leq \lambda_i(x)\leq \beta\leq \lambda_+\quad\text{for any }i=1,\ldots,N\\
\displaystyle{\sum_{i=1}^N\frac{1}{\lambda_i(x)-\alpha_1}\leq \frac{N}{\alpha-\alpha_1}\leq \frac{1}{\lambda_--\alpha_1}\leq \frac{1}{\lambda_--\alpha_1} +\frac{N-1}{\lambda_+-\alpha_1}}\\
\displaystyle{\sum_{i=1}^N\frac{1}{\beta_1-\lambda_i(x)}\leq \frac{N}{\beta_1-\beta}\leq \frac{N-1}{\beta_1-\lambda_+}\leq \frac{1}{\beta_1-\lambda_-} +\frac{N-1}{\beta_1-\lambda_+}}\end{array}\right.
\end{equation}
where 
 $0<\lambda_1(x)\leq\ldots\leq \lambda_N(x)$ are the eigenvalues of $A(x)$ and
$$\lambda_+=\theta\alpha_1+(1-\theta)\beta_1\quad\text{and}\quad
\lambda_-=\displaystyle{\left(\frac{\theta}{\alpha_1}+\frac{1-\theta}{\beta_1}\right)^{-1}}.
$$

A simple computation shows that this is true by taking $\theta=1/3$ for any $N\geq 2$ and $\varepsilon$ small enough.\cvd

\bigskip

We conclude this brief review of properties of $H$- and $G$-convergence by the classical result that $L^1_{loc}$ convergence implies $H$-convergence. At the end of Appendix~\ref{Lpsection} we shall provide for the convenience of the reader a simple proof of this result.

\begin{prop}\label{L1impliesH}
Let $\Omega$ be a bounded open set in $\mathbb{R}^N$, $N\geq 2$.

Let us consider a sequence of conductivity tensors $\{A_n\}_{n\in\mathbb{N}}\subset \tilde{\mathcal{M}}(\alpha,\tilde{\beta})$ and a conductivity tensor $A$ in the same set.

As $n\to\infty$,
if $A_n$ converges to $A$ in $L^1_{loc}(\Omega)$ then $A_n$ $H$-converges to $A$.
\end{prop}

\section{Variational methods of reconstruction}
\label{applicationsection}

Throughout this section $\Omega\subset \mathbb{R}^N$, $N\geq 2$, will be a fixed bounded domain with Lipschitz boundary.

Let us assume that we have a measured, therefore approximated, Dirichlet-to-Neumann map 
$\hat{\Lambda}_D\in  \mathcal{L}(W^{1/2,2}(\partial \Omega),W^{-1/2,2}_{\ast}(\partial \Omega))$
 or Neumann-to-Dirichlet map $\hat{\Lambda}_N\in  \mathcal{L}(W^{-1/2,2}_{\ast}(\partial \Omega),W^{1/2,2}_{\ast}(\partial \Omega))$.

Then a variational approach to the inverse conductivity problem is to solve one of the following minimization problems
$$\min_{A\in \mathcal{M}}\|D\text{-}N(A)-\hat{\Lambda}_D\|\quad\text{or}\quad
\min_{A\in \mathcal{M}}\|N\text{-}D(A)-\hat{\Lambda}_N\|$$
where $\mathcal{M}$ is a suitable class of conductivity tensors and suitable norms are to be chosen.
Alternatively, one may consider the following
$$\min_{A\in \mathcal{M}}\|KD_i(A)\|\quad\text{or}\quad
\min_{A\in \mathcal{M}}\|KN_i(A)\|$$
for $i=1,2$, where, for $i=1$, $\mathcal{M}$ is a suitable class of symmetric conductivity tensors.

Our aim is to show in which cases these minimization problems admit a solution, allowing the class $\mathcal{M}$ to be as large as possible. We shall apply the direct method. Compactness will be provided by the compactness properties of $H$-convergence. Here we shall state corresponding lower semicontinuity results and therefore existence results for our minimization problems. Most proofs will be postponed to the next section.
We begin with the following lower semicontinuity result.

\begin{teo}\label{semicontinuitycor}

Let us consider a sequence of conductivity tensors $\{A_n\}_{n\in\mathbb{N}}\subset \tilde{\mathcal{M}}(\alpha,\tilde{\beta})$ and a conductivity tensor $A$ in the same set.
Let us assume that $A_n$ $H$-converges to $A$ as $n\to\infty$.
%
Then we have
$$\|D\text{-}N(A)-\hat{\Lambda}_D\|
\leq\liminf_n \|D\text{-}N(A_n)-\hat{\Lambda}_D\|$$
where $\|\cdot\|=\|\cdot\|_{\mathcal{L}(W^{1/2,2}(\partial \Omega),W^{-1/2,2}_{\ast}(\partial \Omega))}$. Moreover, for any $i=1,2$,
$$\|KD_i(A)\|
\leq\liminf_n \|KD_i(A_n)\|$$
where $\|\cdot\|=\|\cdot\|_{\mathcal{L}(W^{1/2,2}(\partial \Omega),L^2(\Omega,\mathbb{R}^N))}$.

We also have
$$\|N\text{-}D(A)-\hat{\Lambda}_N\|\leq\liminf_n \|N\text{-}D(A_n)-\hat{\Lambda}_N\|$$
where $\|\cdot\|=\|\cdot\|_{\mathcal{L}(W^{-1/2,2}_{\ast}(\partial \Omega),W^{1/2,2}_{\ast}(\partial \Omega))}$. Moreover, for any $i=1,2$,
$$\|KN_i(A)\|
\leq\liminf_n \|KN_i(A_n)\|$$
where $\|\cdot\|=\|\cdot\|_{\mathcal{L}(W^{-1/2,2}_{\ast}(\partial \Omega),L^2(\Omega,\mathbb{R}^N))}$.

Recall that, both for the Dirichlet and Neumann cases, for $i=1$ we always consider $A_n$ symmetric for any $n\in\mathbb{N}$ and consequently $A$ symmetric as well.
\end{teo}

Actually, by modifying the norms used, 
in certain occasions continuity and not only lower semicontinuity holds.

\begin{teo}\label{Hcontinuity}
Under the assumptions of Theorem~\textnormal{\ref{semicontinuitycor}},
we have
%
$$\lim_{n}\| N\text{-}D(A_n)-N\text{-}D(A)\|_{\mathcal{L}(L^2_{\ast}(\partial\Omega),L^2_{\ast}(\partial\Omega))}=0$$
and
$$\lim_n\|KN_1(A_n)\|_{\mathcal{L}(L^2_{\ast}(\partial\Omega),L^2(\Omega,\mathbb{R}^N))}=\|KN_1(A)\|_{\mathcal{L}(L^2_{\ast}(\partial\Omega),L^2(\Omega,\mathbb{R}^N))}.
$$

\end{teo}

Theorem~\ref{semicontinuitycor} is crucial in establishing the following existence results for the minimum problems stated above

\begin{teo}\label{existence} 
Let us consider the class of conductivity tensors
$\tilde{\mathcal{M}}(\alpha,\tilde{\beta})$. Then there exists a solution to the following minimum problems.

For the Dirichlet-to-Neumann case,
$$\min_{A\in\tilde{\mathcal{M}}(\alpha,\tilde{\beta})}
\|D\text{-}N(A)-\hat{\Lambda}_D\|
$$
where $\|\cdot\|=\|\cdot\|_{\mathcal{L}(W^{1/2,2}(\partial \Omega),W^{-1/2,2}_{\ast}(\partial \Omega))}$, and,
for any $i=1,2$,
$$\min_{A\in\tilde{\mathcal{M}}(\alpha,\tilde{\beta})}
\|KD_i(A)\|
$$
where $\|\cdot\|=\|\cdot\|_{\mathcal{L}(W^{1/2,2}(\partial \Omega),L^2(\Omega,\mathbb{R}^N))}$
and $\tilde{\mathcal{M}}(\alpha,\tilde{\beta})$ has to be replaced by $\mathcal{M}_{sym}(\alpha,\beta)$ for $i=1$.

For the Neumann-to-Dirichlet case,
$$\min_{A\in\tilde{\mathcal{M}}(\alpha,\tilde{\beta})}\|N\text{-}D(A)-\hat{\Lambda}_N\|$$
where $\|\cdot\|=\|\cdot\|_{\mathcal{L}(W^{-1/2,2}_{\ast}(\partial \Omega),W^{1/2,2}_{\ast}(\partial \Omega))}$, and,
for any $i=1,2$,
$$\min_{A\in\tilde{\mathcal{M}}(\alpha,\tilde{\beta})}
\|KN_i(A)\|
$$
where $\|\cdot\|=\|\cdot\|_{\mathcal{L}(W^{-1/2,2}_{\ast}(\partial \Omega),L^2(\Omega,\mathbb{R}^N))}$ and $\tilde{\mathcal{M}}(\alpha,\tilde{\beta})$ has to be replaced by $\mathcal{M}_{sym}(\alpha,\beta)$ for $i=1$.
\end{teo}

\proof{.} Immediate by the compactness properties of $H$-convergence and the previous lower semicontinuity result.\cvd

\bigskip

In the previous theorem we may replace the class of conductivity tensors
$\tilde{\mathcal{M}}(\alpha,\tilde{\beta})$ with the class of symmetric conductivity tensors
$\mathcal{M}_{sym}(\alpha,\beta)$, with exactly the same proof. As we shall see, 
uniqueness is never achieved.

Unfortunately, in general we can not replace $\mathcal{M}_{sym}(\alpha,\beta)$ with
$\mathcal{M}_{scal}(\alpha,\beta)$ and still get existence. We begin by showing that it is not possible to use the direct method of calculus of variations in the case of scalar conductivities.
In fact, let us consider the following argument, limiting ourselves to Neumann-to-Dirichlet maps, since all the other cases are completely analogous. Let us a take a minimizing sequence $\sigma_nI_N\in\mathcal{M}_{scal}(\alpha,\beta)$.Without loss of generality we may assume that $\sigma_nI_N$ $H$-converges, as $n\to \infty$, to $\hat{A}\in \mathcal{M}_{sym}(\alpha,\beta)$. We can not guarantee that $\hat{A}$ is a scalar conductivity, therefore we can not guarantee that we have a minimizer in the class of scalar conductivities unless there exists a scalar conductivity $\sigma I_N$ such that $N\text{-}D(\sigma I_N)=
N\text{-}D(\hat{A})$. We shall show in Example~\ref{Giovanni} that for certain conductivity tensors $\hat{A}$
this may not happen. Here instead we state a nonexistence result for our minimization problems in the class of scalar conductivities.

\begin{exam}\label{nonexistence}
Let us consider $\Omega=B_1\subset \mathbb{R}^2$ and two positive constants $a$ and $b$, with $a\neq b$. Let
\begin{equation}\label{Giovannicond}
\hat{A}(x)=\hat{A}=\left[\begin{matrix}a& 0\\0 &b\end{matrix}\right]\quad\text{for any  }x\in B_1.
\end{equation}
We assume that for some positive constants $\alpha$ and $\beta$ we have
$\alpha\leq a,b\leq\beta$, hence
$\hat{A}\in \mathcal{M}_{sym}(\alpha,\beta)\backslash \mathcal{M}_{scal}(\alpha,\beta)$.

Let us set $\hat{\Lambda}_N=N\text{-}D(\hat{A})$. Let $\alpha_1$ and $\beta_1$ as defined in Proposition~\ref{approxbyscalar}. Then there is no solution to the following minimum problems
$$\min_{A\in\mathcal{M}_{scal}(\alpha_1,\beta_1)}
\|N\text{-}D(A)-\hat{\Lambda}_N\|,
$$
where this time we set $\|\cdot\|=\|\cdot\|_{
_{\mathcal{L}(L^2_{\ast}(\partial\Omega),L^2_{\ast}(\partial\Omega))}}$,
and
$$\min_{A\in\mathcal{M}_{scal}(\alpha_1,\beta_1)}
\|KN_1(A)\|,
$$
where this time we set $\|\cdot\|=\|\cdot\|_{
_{\mathcal{L}(L^2_{\ast}(\partial\Omega),L^2(\Omega,\mathbb{R}^N))}}$.
\end{exam}

\subsubsection{Finite number of measurements}

In the sequel of this section for simplicity we limit ourselves to symmetric conductivity tensors. 

A variety of variational methods have been proposed in the literature for solving the inverse conductivity problem. We limit ourselves to two classical approaches, namely the one of Yorkey, \cite{Yor}, and the ones by Kohn and Vogelius, \cite{Koh e Vog87}, and Kohn and McKenney, \cite{Koh e McK}.

All of them make use of a finite number of measurements, therefore let us consider the following
common setting. We assume that $n$ different measurements are performed, namely we prescribe $g_1,\ldots g_n$ different current densities at the boundary and we measure the corresponding potentials $\varphi_1,\ldots \varphi_n$ again at the boundary. This is the more commonly accepted framework 
for performing electrostatic measurements on the boundary,
even if a more precise approach from an application point of view is the so-called experimental measurements model which we shall discuss next in this section.

It is not restrictive to assume that
all the current densities $g_i$, $i=1,\ldots,n$, belongs to $L^2_{\ast}(\partial\Omega)$.

Then the approach of Yorkey is to consider the following minimization problem,
$$\min_{A\in \mathcal{M}}\sum_{i=1}^n\int_{\partial\Omega}|N\text{-}D(A)(g_i)-\varphi_i |^2=\min_{A\in \mathcal{M}}J_0(A)$$
where $\mathcal{M}$ is a suitable class of conductivity tensors.

Kohn, Vogelius and McKenney, on the other hand, proposed the following minimization problems. Either
$$\min_{A\in \mathcal{M}}\sum_{i=1}^n\int_{\Omega}|KN_1(A)[g_i]|^2=\min_{A\in \mathcal{M}}J_1(A)$$
or 
$$\min_{A\in \mathcal{M}}\sum_{i=1}^n\int_{\Omega}|KN_2(A)[g_i]|^2=\min_{A\in \mathcal{M}}J_2(A)$$
where we assume that $\hat{\Lambda}_N[g_i]=\varphi_i\in W^{1/2,2}_{\ast}(\partial\Omega)$.
We notice that this approach has the disadvantage that the measured potentials $\varphi_i$, $i=1,\ldots,n$, must belong to $W^{1/2,2}_{\ast}(\partial\Omega)$, whereas in the approach of Yorkey it is enough that they belong to $L^2_{\ast}(\partial\Omega)$. On the other hand, assuming that for any $i$, $\varphi_i$ is an interpolation of pointwise values of the potential at a finite number of fixed points on the boundary, and assuming the boundary smooth enough, a piecewise linear interpolation would do for both cases.

By the previous results in this section it is easy to show that $J_0$, $J_1$ and $J_2$ are lower semicontinuous on $\mathcal{M}=\mathcal{M}_{sym}(\alpha,\beta)$ with respect to the $G$-convergence, by Theorem~\ref{semicontinuitycor}. Moreover,
since we are assuming that all the current densities $g_i$, $i=1,\ldots,n$, belong to $L^2_{\ast}(\partial\Omega)$,
we infer by Theorem~\ref{Hcontinuity}  that $J_0$ and $J_1$ are actually continuous.
By Proposition~\ref{L1impliesH} the same semicontinuity or continuity properties hold if we replace the $G$-convergence by the $L^1$ convergence.
By compactness of $\mathcal{M}_{sym}(\alpha,\beta)$ with respect to the $G$-convergence, we have existence of a minimizer, for any of these three $J_i$, $i=0,1,2$, functionals, if we set
$\mathcal{M}=\mathcal{M}_{sym}(\alpha,\beta)$. On the other hand, we might expect nonexistence issues, as previously shown, if we set the minimization problems
on the class $\mathcal{M}=\mathcal{M}_{scal}(\alpha,\beta)$. However using symmetric conductivity tensors instead of scalar ones in the numerical reconstruction has several drawbacks, as discussed in \cite{Koh e McK}.

\subsubsection{Experimental measurements}

Similar continuity properties, leading to existence results for corresponding minimum problems, hold also for the so-called \emph{experimental measurements}, which are
measurements which can be actually obtained from the experiments
and have been modeled
in \cite{Som e Che e Isa}.

We describe the model, for further details we refer to the original paper. We assume that the conductor is contained in a bounded domain $\Omega\subset\mathbb{R}^N$, $N\geq 2$, with Lipschitz boundary, and that $A$ is a symmetric conductivity tensor in $\Omega$, namely $A\in \mathcal{M}_{sym}(\alpha,\beta)$.
We attach $L$ \emph{electrodes} on the boundary of the conductor $\Omega$.
The \emph{contact regions} between the electrodes and the conductor are subsets of
$\partial \Omega$ and will
be denoted
by $e_l$, $l=1,\ldots,L$. We assume that the subsets $e_l$, $l=1,\ldots,L$,
are open, nonempty, connected, with a smooth boundary and such that
their closures are pairwise disjoint.
Any electrode is identified with its contact region.
A current is sent to the body through the electrodes and the corresponding voltages
are measured on the same electrodes. For each $l$, $l=1,\ldots,L$,
the current applied to the electrode $e_l$ will be denoted by
$I_l$ and the voltage measured on the electrode will be denoted by $V_l$. The column vector $I$ whose
components are $I_l$, $l=1,\ldots,L$, is a \emph{current pattern}
if the condition $\sum_{l=1}^LI_l=0$ is satisfied.
The corresponding \emph{voltage pattern}, that is the column vector $V$ whose components are $V_l$, $l=1,\ldots,L$,
is determined up to an additive constant and we always choose to normalize it in such a way that
$\sum_{l=1}^LV_l=0$. The voltage pattern depends on the current pattern in a linear way, through an
$L\times L$ symmetric matrix $R=R(A)$ which is called the \emph{resistance matrix}, that is
$V=RI$. Without loss of generality we assume that 
$R[1]=0$, where $[1]$ denotes the column vector whose components are all equal to $1$.

We briefly recall the model
used to determine the resistance matrix $R$. We assume that at each
electrode $e_l$, $l=1,\ldots,L$, a \emph{surface impedance} is present and we denote it with $z_l$.
Let us assume that there exists $Z_1$, $Z_2$, $0<Z_1<Z_2$, such that
for each $l$, $l=1,\ldots,L$, $Z_1\leq z_l\leq Z_2$.
If we apply the current pattern
$I$ on the electrodes, then the voltage $u$ inside the body satisfies the following
boundary value problem
\begin{equation}\label{expmeaspbm}
\left\{\begin{array}{ll}
\mathrm{div}(A\nabla u)=0 &\text{in }\Omega,\\
u+z_l A\nabla u\cdot\nu=U_l &\text{on }e_l,\ l=1,\ldots,L,\\
A\nabla u\cdot\nu=0 &\text{on }\partial\Omega\backslash\bigcup_{l=1}^Le_l,\\
\int_{e_l}A\nabla u\cdot\nu=I_l&\text{for any }l=1,\ldots,L,
\end{array}\right.
\end{equation}
where $U_l$, $l=1,\ldots,L,$ are constants to be determined. We call $U$ the column vector whose
components are given by $U_l$, $l=1,\ldots,L$.

For any $l$, $l=1,\ldots,L$, $V_l$, a component of the voltage pattern $V$, is given by
$V_l=\int_{e_l}u$, thus, by \eqref{expmeaspbm},
$$V_l=\mathcal{H}^{N-1}(e_l)U_l-z_lI_l,$$
where $\mathcal{H}^{N-1}$ denotes the $(N-1)$-dimensional Hausdorff measure.

Let $H=W^{1,2}(\Omega)\times \mathbb{R}^L$.
By \cite[Theorem~3.3]{Som e Che e Isa}, there exists a unique couple $(u,U)\in H$,
satisfying
$\sum_{l=1}^L\mathcal{H}^{N-1}(e_l)U_l-z_lI_l=0$, such that \eqref{expmeaspbm} is satisfied. Thus the current pattern
$I$ uniquely determines the voltage pattern $V$, if this is normalized in such a way that
$\sum_{l=1}^LV_l=0$. Furthermore, it has been proved in \cite{Som e Che e Isa} that the relation
between $I$ and $V$ is linear, thus the resistance matrix $R(A)$ is well defined, and that $R(A)$ is actually symmetric.

We briefly recall the argument. For any symmetric conductivity tensor $A$ and any
$(u,U)$, $(w,W)\in H$, we let 
$$B_{A}((u,U),(w,W))=\int_{\Omega}A\nabla u\cdot\nabla w+
\sum_{l=1}^L\frac{1}{z_l}\int_{e_l}(u-U_l)(w-W_l).$$
Then, for any $g\in W^{-1/2,2}_{\ast}(\partial\Omega)$ and any $I\in\mathbb{R}^L$
such that $\sum_{l=1}^LI_l=0$, we have that there exists a couple $(u,U)\in H$ satisfying
\begin{equation}\label{eqform}
B_{A}((u,U),(w,W))=\langle g, w\rangle+\sum_{l=1}^LI_lW_l\quad\text{for any }(w,W)\in H.
\end{equation}
The solution is unique up to adding the same constant to both $u$ and $U$.

We can find a constant $C_1$, depending on $N$, $\Omega$, $\alpha$, $Z_2$ and the electrodes only,
such that
\begin{equation}\label{coerc}
\|\nabla u\|_{L^2(\Omega)}\leq C_1(\|g\|_{W^{-1/2,2}_{\ast}(\partial\Omega)}+\|I\|).
\end{equation}
Notice that if $(u,U)$ solves \eqref{eqform}, then
\begin{equation}\label{diveq}
\mathrm{div}(A\nabla u)=0\quad\text{in }\Omega.
\end{equation}
We 
denote
$\phi=A\nabla u\cdot\nu\in W^{-1/2,2}_{\ast}(\partial\Omega)$
and recall that
$N\text{-}D(A)$ is the Neumann-to-Dirichlet map associated to $A$.
Then $(u,U)$
solves \eqref{eqform} if and only if $u$ satisfies \eqref{diveq}, we have
\begin{equation}\label{Uexp}
\mathcal{H}^{N-1}(e_l)U_l=z_lI_l+
\int_{e_l}u,
\quad\text{for any }l=1,\ldots,L,
\end{equation}
and the following equation holds in $W^{-1/2,2}_{\ast}(\partial\Omega)$
\begin{multline}\label{boundaryeqform}
\phi+\sum_{l=1}^L\frac{1}{z_l}\left(N\text{-}D(A)[\phi]-
\frac{1}{\mathcal{H}^{N-1}(e_l)}\int_{e_l}N\text{-}D(A)[\phi]\right)\chi_{e_l}
=\\ g+\sum_{l=1}^L\left(\frac{I_l}{\mathcal{H}^{N-1}(e_l)}\chi_{e_l}\right).
\end{multline}

Let $\mathcal{K}(A):W^{-1/2.2}_{\ast}(\partial\Omega)\to L^2_{\ast}(\partial\Omega)$
be the
operator defined as follows. For any $\phi\in W^{-1/2,2}_{\ast}(\partial\Omega)$
$$\mathcal{K}(A)[\phi]=\sum_{l=1}^L\frac{1}{z_l}\left(N\text{-}D(A)[\phi]-
\frac{1}{\mathcal{H}^{N-1}(e_l)}\int_{e_l}N\text{-}D(A)[\phi]\right)\chi_{e_l}.$$
We have that $\mathcal{K}(A)$ is a compact linear operator, and, for a constant 
$C_2$ depending on $N$, $\Omega$, $\alpha$, $Z_1$ and the electrodes only, we have
\begin{equation}\label{normk}
\|\mathcal{K}(A) \|_{\mathcal{L}(W^{-1/2.2}_{\ast}(\partial\Omega),L^2_{\ast}(\partial\Omega))}\leq C_2.
\end{equation} 
Therefore,
$\mathcal{K}(A)$ is a compact linear operator also from $W^{-1/2.2}_{\ast}(\partial\Omega)$
into itself and from $L^2_{\ast}(\partial\Omega)$ into itself.

Since for any $g\in W^{-1/2.2}_{\ast}(\partial\Omega)$ and $I=0$
the equation \eqref{eqform} admits a solution, we can infer that
$Id+\mathcal{K}(A):W^{-1/2.2}_{\ast}(\partial\Omega)\to W^{-1/2.2}_{\ast}(\partial\Omega)$
is bijective, $Id$ denoting the identity operator.
We deduce that 
$Id+\mathcal{K}(A):L^2_{\ast}(\partial\Omega)\to L^2_{\ast}(\partial\Omega)$ is bijective as well.
We denote with $\tilde{K}(A)$ the inverse to $Id+\mathcal{K}(A)$. Using \eqref{coerc}
we have\begin{equation}\label{unifest2}
\|\tilde{\mathcal{K}}(A)[g]\|_{W^{-1/2.2}_{\ast}(\partial\Omega)}\leq C_3\|g\|_{W^{-1/2.2}_{\ast}(\partial\Omega)},
\end{equation}
where $C_3$ depends
on $N$, $\Omega$, $\alpha$, $\beta$, $Z_2$ and the electrodes only.

We conclude that, if $g\in L^2_{\ast}(\partial\Omega)$, then, for a constant $C_4$ depending on 
$N$, $\Omega$, $\alpha$, $\beta$, $Z_2$ and the electrodes only, we have
$$\|\tilde{\mathcal{K}}(A)[g]\|_{W^{-1/2.2}_{\ast}(\partial\Omega)}\leq C_4\|g\|_{L^2_{\ast}(\partial\Omega)},
$$
hence
\begin{multline}\label{L2L2}
\|\tilde{\mathcal{K}}(A)[g]\|_{L^2_{\ast}(\partial\Omega)}\leq\\
\|\mathcal{K}(A) [|\tilde{\mathcal{K}}(A)[g]] \|_{L^2_{\ast}(\partial\Omega)}+
\|g\|_{L^2_{\ast}(\partial\Omega)}
\leq (C_2 C_4 +1)\|g\|_{L^2_{\ast}(\partial\Omega)}.
\end{multline}

For any given current pattern $I$, that is $I\in \mathbb{R}^L$ such that
$\sum_{l=1}^LI_l=0$, we can define $\tilde{I}=\sum_{l=1}^L\left(\frac{I_l}{\mathcal{H}^{N-1}(e_l)}\chi_{e_l}\right)\in
L^2_{\ast}(\partial\Omega)$. Furthermore, there exists a constant $C_5$, depending on $N$, $\Omega$ and
the electrodes only, such that
\begin{equation}\label{ItildeI}
\|\tilde{I}\|_{L^2_{\ast}(\partial\Omega)}\leq C_5\|I\|.
\end{equation}
As pointed out in \cite{Som e Che e Isa},
we have that $(u,U)$ solves our direct problem \eqref{expmeaspbm}
for a given current pattern $I$
if and only if
\eqref{eqform} is satisfied with $g=0$. Therefore, we have that $R(A)I=V$ where, for any $l=1,\ldots,L$,  
\begin{equation}\label{R-Kconn}
V_l=\int_{e_l}N\text{-}D(A)(\tilde{\mathcal{K}}(A)\tilde{I})+c\mathcal{H}^{N-1}(e_l),
\end{equation}
where $c$ is a constant which can be computed by imposing the condition
that $\sum_{l=1}^LV_l=0$, that is
\begin{equation}\label{constantc}
c=-\frac{\sum_{l=1}^L\int_{e_l}N\text{-}D(A)(\tilde{\mathcal{K}}(A)\tilde{I})}{\sum_{l=1}^L\mathcal{H}^{N-1}(e_l)}.
\end{equation}

Hence we can define our forward operator $R:\mathcal{M}_{sym}(\alpha,\beta)\to\mathbb{M}^{N\times N}(\mathbb{R})$ such that to any $A\in \mathcal{M}_{sym}(\alpha,\beta)$ it associates its corresponding resistance matrix $R(A)$.

The following result holds.

\begin{prop}\label{expmeasprop}
Let $A_1$, $A_2\in \mathcal{M}_{sym}(\alpha,\beta)$.
Then there exists a constant $C$,
depending on 
$N$, $\Omega$, $\alpha$, $\beta$, $Z_1$, $Z_2$ and the electrodes only,
such that
\begin{equation}\label{expmeasest}
\|R(A_1)-R(A_2)\|\leq C\|N\text{-}D(A_1)-N\text{-}D(A_2)\|_{\mathcal{L}(L^2_{\ast}(\partial\Omega),L^2_{\ast}(\partial\Omega))}.
\end{equation}
\end{prop}

\proof{.}
For any $I\in\mathbb{R}^L$ such that
$\sum_{l=1}^LI_l=0$ we evaluate
$\|(R(A_1)-R(A_2))I\|$. We recall that we have set $R(A_1)[1]=R(A_2)[1]=0$.
By \eqref{R-Kconn} and \eqref{constantc},
we have that
$$\|(R(A_1)-R(A_2))I\|\leq C_6\|N\text{-}D(A_1)(\tilde{\mathcal{K}}(A_1)\tilde{I})-
N\text{-}D(A_2)(\tilde{\mathcal{K}}(A_2)\tilde{I})\|_{L^2_{\ast}(\partial\Omega)},$$
where $C_6$ depends on $N$ and the electrodes only. Thus
\begin{multline*}
\|(R(A_1)-R(A_2))I\| \leq C_6\big(\|N\text{-}D(A_1)(\tilde{\mathcal{K}}(A_1)-\tilde{\mathcal{K}}(A_2))
\tilde{I}\|_{L^2_{\ast}(\partial\Omega)}+\\
\|(N\text{-}D(A_1)-N\text{-}D(A_2))(\tilde{\mathcal{K}}(A_2)
\tilde{I})\|_{L^2_{\ast}(\partial\Omega)}\big),
\end{multline*}
and, by \eqref{L2L2} and \eqref{ItildeI}, we obtain that
\begin{multline}\label{Rest}
\|(R(A_1)-R(A_2))I\| \leq
C_6C_5\tilde{C}\|\tilde{\mathcal{K}}(A_1)-\tilde{\mathcal{K}}(A_2)\|_{\mathcal{L}(L^2_{\ast}(\partial\Omega),L^2_{\ast}(\partial\Omega))}
\|I\|+\\
 C_6C_5\tilde{C}_1\|N\text{-}D(A_1)-N\text{-}D(A_2)\|_{\mathcal{L}(L^2_{\ast}(\partial\Omega),L^2_{\ast}(\partial\Omega))}\|I\|.
\end{multline}
where $\tilde{C}$ depends on $N$, $\Omega$ and $\alpha$ only and $\tilde{C}_1=C_2 C_4 +1$

It remains to evaluate the term
$\|\tilde{\mathcal{K}}(A_1)-\tilde{\mathcal{K}}(A_2)\|_{\mathcal{L}(L^2_{\ast}(\partial\Omega),L^2_{\ast}(\partial\Omega))}$. We proceed as follows.
Using the identity
$$\tilde{\mathcal{K}}(A_1)-\tilde{\mathcal{K}}(A_2)=\tilde{\mathcal{K}}(A_1)[\mathcal{K}(A_2)-\mathcal{K}(A_1)]\tilde{K}(A_2)
$$
we obtain that
\begin{equation}\label{K1-K2}
\|\tilde{\mathcal{K}}(A_1)-\tilde{\mathcal{K}}(A_2)\|_{\mathcal{L}(L^2_{\ast}(\partial\Omega),L^2_{\ast}(\partial\Omega))}\leq \tilde{C}_1^2
\|\mathcal{K}(A_1)-\mathcal{K}(A_2)\|_{\mathcal{L}(L^2_{\ast}(\partial\Omega),L^2_{\ast}(\partial\Omega))}.
\end{equation}
It is not difficult to show that for some constant $\tilde{C}_2$, depending on
$N$, $\Omega$, $Z_1$ and the electrodes only,
we have
$$\|\mathcal{K}(A_1)-\mathcal{K}(A_2)\|_{\mathcal{L}(L^2_{\ast}(\partial\Omega),L^2_{\ast}(\partial\Omega))}\leq \tilde{C}_2
\|N\text{-}D(A_1)-N\text{-}D(A_2)\|_{\mathcal{L}(L^2_{\ast}(\partial\Omega),L^2_{\ast}(\partial\Omega))}
$$
therefore the proof is concluded.\cvd

\bigskip

As an immediate corollary of Proposition~\ref{expmeasprop}, we obtain that, on $\mathcal{M}=\mathcal{M}_{sym}(\alpha,\beta)$, the forward operator $R$ is continuous with respect to the $G$-conver\-gence, by Theorem~\ref{Hcontinuity}. Moreover, $R$ is
indeed H\"older continuous with respect to the $L^1$ norm on $\mathcal{M}$, by Theorem~\ref{pcontinuity} which we shall state in the Appendix.

\section{Continuity with respect to $H$-convergence}\label{Hsection}

Let $\Omega\subset \mathbb{R}^N$, $N\geq 2$, be a bounded open set. As before, for any $\varphi\in W^{1,2}(\Omega)$, $F\in (W^{1,2}\Omega))^{\ast}$ and any conductivity tensor $A$, we let $u=u(A,F,\varphi)$ be the unique weak solution to \eqref{Dirichletproblem}.

Assuming further that $\Omega$ is a domain and it has a Lipschitz boundary, and that $F\in (W^{1,2}(\Omega))^{\ast}$ is such that $\langle F,1\rangle_{((W^{1,2}(\Omega))^{\ast},W^{1,2}(\Omega))}=0$, we let $v=v(A,F)$ be the unique weak solution to \eqref{Neumannproblem}

We begin with the following lemma, that extends, under $H$-convergence of the conductivity tensors, convergence of solutions  to the boundary value problems we are interested in.

\begin{lem}\label{bvpHconv}
Let us consider a sequence of conductivity tensors $\{A_n\}_{n\in\mathbb{N}}\subset \tilde{\mathcal{M}}(\alpha,\tilde{\beta})$ and a conductivity tensor $A$ in the same set. 
Assume that $A_n$ $H$-converges to $A$ as $n\to \infty$.
We denote $u_n=u(A_n,F,\varphi)$ and $u=u(A,F,\varphi)$, and $v_n=v(A_n,F)$ and $v=v(A,F)$

Then we have that
$u_n$ converges to $u$ weakly in $W^{1,2}(\Omega)$ and
$$A_n\nabla u_n\rightharpoonup A\nabla u \quad \text{weakly in }L^2(\Omega,\mathbb{R}^N).$$

Analogously, we have tha $v_n$ converges to $v$ weakly in $W^{1,2}(\Omega)$ and
$$A_n\nabla v_n\rightharpoonup A\nabla v \quad \text{weakly in }L^2(\Omega,\mathbb{R}^N).$$

\end{lem}

\proof{.}
We sketch the proof for the convenience of the reader. We have that, passing to a subsequence, $u_{n_k}$ converges to a function $\tilde{u}$ weakly in $W^{1,2}(\Omega)$ and $A_{n_k}\nabla u_{n_k}$ converges to a vector $V$ weakly in $L^2(\Omega,\mathbb{R}^N)$. Clearly $\tilde{u}=\varphi$ on $\partial\Omega$ in a weak sense. By the independence of the $H$-convergence from the boundary data, \cite[Theorem~1]{Mur e Tar1}, we have that
$$A_{n_k}\nabla u_{n_k}\rightharpoonup A\nabla \tilde{u}$$ weakly in $L^2(\Omega_1,\mathbb{R}^N)$, for any open $\Omega_1$ compactly contained in $\Omega$. Therefore $V=A\nabla \tilde{u}$
and we may conclude that $\tilde{u}=u(A,F,\varphi)$. It is now easy to conclude that the whole sequences $u_n$ and $A_n\nabla u_n$ converge to $u$ and $A\nabla u$, respectively, as required.
We conclude that, for any $\psi\in W^{1,2}(\Omega)$ we have, as $n\to \infty$,
$$\int_{\Omega}A_n\nabla u_n\cdot\nabla \psi\to \int_{\Omega}A\nabla u\cdot\nabla \psi.$$

For Neumann problems, the proof is completely analogous.\cvd

\bigskip

We easily apply these previous results to Dirichlet-to-Neumann maps or Neumann-to-Dirichlet maps and obtain the following continuity properties.

\begin{prop}\label{continuityprop}
Let $\Omega\subset \mathbb{R}^N$, $N\geq 2$, be a bounded domain with Lipschitz boundary.
We consider a sequence of conductivity tensors $\{A_n\}_{n\in\mathbb{N}}\subset \tilde{\mathcal{M}}(\alpha,\tilde{\beta})$ and a conductivity tensor $A$ in the same set.

%

If $A_n$ $H$-converges to $A$ as $n\to\infty$, then
for any $\varphi$, $\psi\in W^{1/2,2}(\partial\Omega)$ we have
$$(D\text{-}N(A_n)(\varphi)-D\text{-}N(A)(\varphi))[\psi]\to 0$$
and for any $g\in W^{-1/2,2}_{\ast}(\partial\Omega)$ we have
$$(N\text{-}D(A_n)(g)-N\text{-}D(A)(g))\rightharpoonup 0\quad\text{weakly in }W_{\ast}^{1/2,2}(\partial\Omega).$$

Finally,
for any $\varphi\in W^{1/2,2}(\partial\Omega)$ we have
\begin{equation}\label{encons1}
\int_{\Omega}A_n\nabla u_n\cdot\nabla u_n\to \int_{\Omega}A\nabla u\cdot\nabla u
\end{equation}
where $u_n=u(A_n,0,\varphi)$ and $u=u(A,0,\varphi)$, solutions to \eqref{Dirichletproblem1}, and
for any $g\in W^{-1/2,2}_{\ast}(\partial\Omega)$ we have
\begin{equation}\label{encons2}
\int_{\Omega}A_n\nabla v_n\cdot\nabla v_n\to \int_{\Omega}A\nabla v\cdot\nabla v
\end{equation}
where $v_n=v(A_n,g)$ and $v=v(A,g)$, solutions to \eqref{Neumannproblem1}.
\end{prop}

As a corollary we may prove Theorem~\ref{semicontinuitycor} and Theorem~\ref{Hcontinuity}.

\proof{ of Theorem~\ref{semicontinuitycor}.} It follows by the Banach-Steinhaus Theorem. In fact, for any fixed $\varphi\in  W^{1/2,2}(\partial\Omega)$ let us call $F_n(\varphi)=(D\text{-}N(A_n)(\varphi)-\hat{\Lambda}_D(\varphi))\in W^{-1/2,2}_{\ast}(\partial\Omega)$, $n\in\mathbb{N}$, and $F(\varphi)=(D\text{-}N(A)(\varphi)-\hat{\Lambda}_D(\varphi))\in W^{-1/2,2}_{\ast}(\partial\Omega)$. Since for any $\psi\in  W^{1/2,2}(\partial\Omega)$ we have
$$\lim_n F_n(\varphi)[\psi]=F(\varphi))[\psi]$$
we conclude that
$$\|F(\varphi)\|_{W^{-1/2,2}_{\ast}(\partial\Omega)}\leq\liminf_n \|F_n(\varphi)\|_{W^{-1/2,2}_{\ast}(\partial\Omega)}.$$
Hence for any $\varphi\in  W^{1/2,2}(\partial\Omega)$
\begin{multline*}
\|(D\text{-}N(A)-\hat{\Lambda}_D)(\varphi)\|_{W^{-1/2,2}_{\ast}(\partial\Omega)}\leq\\
\liminf_n
\|(D\text{-}N(A_n)-\hat{\Lambda}_D)\|_{\mathcal{L}(W^{1/2,2}(\partial \Omega),W^{-1/2,2}_{\ast}(\partial \Omega))}
 \|\varphi\|_{W^{1/2,2}(\partial\Omega)}
 \end{multline*}
and the lower semicontinuity of the Dirichlet-to-Neumann forward operator
is proved. Passing to $KD_i$, $i=1,2$, using a similar argument, it is enough to notice that 
$$\|KD_2(A)[\varphi] \|_{L^2(\Omega,\mathbb{R}^N)}\leq \liminf_n \|KD_2(A_n)[\varphi] \|_{L^2(\Omega,\mathbb{R}^N)}$$
and that
$$\|KD_1(A)[\varphi] \|_{L^2(\Omega,\mathbb{R}^N)}=\lim_n \|KD_1(A_n)[\varphi] \|_{L^2(\Omega,\mathbb{R}^N)}.$$
Here we used \eqref{encons1} and \eqref{encons2} and the fact that, setting $u=u(A,0,\varphi)$ and
$v(A,\hat{\Lambda}_D(\varphi))$, we have
\begin{multline*}
\|KD_1(A)[\varphi]\|^2_{L^2(\Omega,\mathbb{R}^N)}=
\int_{\Omega}A\nabla u\cdot\nabla u+
\int_{\Omega}A\nabla v\cdot\nabla v
-2\int_{\Omega}A\nabla u\cdot\nabla v=\\
\int_{\Omega}A\nabla u\cdot\nabla u+
\int_{\Omega}A\nabla v\cdot\nabla v-2\langle \hat{\Lambda}_D(\varphi),\varphi\rangle_{(W^{-1/2,2}(\partial\Omega),W^{1/2,2}(\partial\Omega))}.
\end{multline*}

About the second part, let us fix $g\in W^{-1/2,2}_{\ast}(\partial\Omega)$. Then
$(N\text{-}D(A_n)-\hat{\Lambda}_N)(g)$ converges, as $n\to\infty$, to 
$(N\text{-}D(A)-\hat{\Lambda}_N)(g)$ weakly in $W_{\ast}^{1/2,2}(\partial\Omega)$, therefore
\begin{multline*}
\|(N\text{-}D(A)-\hat{\Lambda}_N)(g)\|_{W_{\ast}^{1/2,2}(\partial\Omega)}\leq
\liminf_n \|(N\text{-}D(A_n)-\hat{\Lambda}_N)(g)\|_{W_{\ast}^{1/2,2}(\partial\Omega)}\leq \\
\liminf_n \|(N\text{-}D(A_n)-\hat{\Lambda}_N)\|_{\mathcal{L}(W^{-1/2,2}_{\ast}(\partial \Omega),W^{1/2,2}_{\ast}(\partial \Omega))}
\|g\|_{W_{\ast}^{-1/2,2}(\partial\Omega)}
\end{multline*}
and the lower semicontinuity of the Neumann-to-Dirichlet forward operator
is proved. About $KN_{i}$ again it is enough to show that
$$\|KN_2(A)[g] \|_{L^2(\Omega,\mathbb{R}^N)}\leq \liminf_n \|KN_2(A_n)[g] \|_{L^2(\Omega,\mathbb{R}^N)}$$
and that
$$\|KN_1(A)[g] \|_{L^2(\Omega,\mathbb{R}^N)}=\lim_n \|KN_1(A_n)[g] \|_{L^2(\Omega,\mathbb{R}^N)},$$
which follows again from \eqref{encons1} and \eqref{encons2} and from the fact
that, setting $u=u(A,0,\hat{\Lambda}_N(g))$ and
$v(A,g)$, we have
\begin{multline*}
\|KN_1(A)[g]\|^2_{L^2(\Omega,\mathbb{R}^N)}=
\int_{\Omega}A\nabla u\cdot\nabla u+
\int_{\Omega}A\nabla v\cdot\nabla v
-2\int_{\Omega}A\nabla u\cdot\nabla v=\\
\int_{\Omega}A\nabla u\cdot\nabla u+
\int_{\Omega}A\nabla v\cdot\nabla v-2\langle g,\hat{\Lambda}_N(g)\rangle_{(W^{-1/2,2}_{\ast}(\partial\Omega),W^{1/2,2}_{\ast}(\partial\Omega))}.
\end{multline*}
The proof is concluded.\cvd

\bigskip

The fact that $W^{1/2,2}_{\ast}(\partial\Omega)\subset L^2_{\ast}(\partial\Omega)$ and 
$L^2_{\ast}(\partial\Omega)\subset W^{-1/2,2}_{\ast}(\partial\Omega)$ with compact immersions, and the density of $L^2_{\ast}(\partial\Omega)$ in $W^{-1/2,2}_{\ast}(\partial\Omega)$, allow us to prove Theorem~\ref{Hcontinuity}.

\proof{ of Theorem~\ref{Hcontinuity}.}
First of all, we notice that there exists a constant $C$ such that
\begin{equation}\label{uniformbound}
\|N\text{-}D(A_n)-N\text{-}D(A)\|_{\mathcal{L}(W^{-1/2,2}_{\ast}(\partial \Omega),W^{1/2,2}_{\ast}(\partial \Omega))}\leq C\quad\text{fo any }n\in\mathbb{N}.
\end{equation}
Furthermore,
for any $g\in W^{-1/2,2}_{\ast}(\partial\Omega)$, we have that
\begin{equation}\label{strongconverg}
\lim_n\|(N\text{-}D(A_n)-N\text{-}D(A))(g)\|_{L^2_{\ast}(\partial\Omega)}=0.
\end{equation}

For any $n\in\mathbb{N}$, let $g_n\in L^2_{\ast}(\partial\Omega)$ be such that
$\|g_n\|_{L^2_{\ast}(\partial\Omega)}=1$ and
\begin{multline*}
\| N\text{-}D(A_n)-N\text{-}D(A)\|_{\mathcal{L}(L^2_{\ast}(\partial\Omega),L^2_{\ast}(\partial\Omega))}-1/n\leq\\ \|(N\text{-}D(A_n)-N\text{-}D(A))(g_n)\|_{L^2_{\ast}(\partial\Omega)}
 \leq \| N\text{-}D(A_n)-N\text{-}D(A)\|_{\mathcal{L}(L^2_{\ast}(\partial\Omega),L^2_{\ast}(\partial\Omega))}.
 \end{multline*}

Let us consider a subsequence such that
\begin{multline*}
\lim_k\| N\text{-}D(A_{n_k})-N\text{-}D(A)\|_{\mathcal{L}(L^2_{\ast}(\partial\Omega),L^2_{\ast}(\partial\Omega))}=\\\limsup_n \| N\text{-}D(A_n)-N\text{-}D(A)\|_{\mathcal{L}(L^2_{\ast}(\partial\Omega),L^2_{\ast}(\partial\Omega))}
\end{multline*}
and such that $g_{n_k}$ converges, as $k\to\infty$, to $g$ strongly in $W^{-1/2,2}_{\ast}(\partial \Omega)$.
Then
\begin{multline*}
\|(N\text{-}D(A_{n_k})-N\text{-}D(A))(g_{n_k})\|_{L^2_{\ast}(\partial\Omega)}
 \leq\\ \|(N\text{-}D(A_{n_k})-N\text{-}D(A))(g_{n_k}-g)\|_{L^2_{\ast}(\partial\Omega)}+
 \|(N\text{-}D(A_{n_k})-N\text{-}D(A))(g)\|_{L^2_{\ast}(\partial\Omega)}.
 \end{multline*}
The second term of the right hand side goes to zero by \eqref{strongconverg}. For what concerns the first term of the right hand side, we notice that, for some constant $C_1$, and using also \eqref{uniformbound}, we have
\begin{multline*}
\|(N\text{-}D(A_{n_k})-N\text{-}D(A))(g_{n_k}-g)\|_{L^2_{\ast}(\partial\Omega)}\leq\\
C_1\|(N\text{-}D(A_{n_k})-N\text{-}D(A))(g_{n_k}-g)\|_{W^{1/2,2}_{\ast}(\partial\Omega)}\leq
C_1C\|g_{n_k}-g\|_{W^{-1/2,2}_{\ast}(\partial\Omega)}
\end{multline*}
and the proof of the first part is concluded.

Let us now consider the second part. Again it is not difficult to show that there exists a constant $C$ such that for any $n\in\mathbb{N}$ we have
$$\| KN_1(A_n)\|_{\mathcal{L}(W^{-1/2,2}_{\ast}(\partial \Omega),L^2(\Omega,\mathbb{R}^N))}
\leq C.$$
We already know that for any $g\in W^{-1/2,2}_{\ast}(\partial\Omega)$, we have that
\begin{equation}\label{strongconverg2}
\|KN_1(A)[g]\|_{L^2(\Omega,\mathbb{R}^N)}=\lim_n\|KN_1(A_n)[g]\|_{L^2(\Omega,\mathbb{R}^N)}.
\end{equation}

For any $n\in\mathbb{N}$, let $g_n\in L^2_{\ast}(\partial\Omega)$ be such that
$\|g_n\|_{L^2_{\ast}(\partial\Omega)}=1$ and
\begin{multline*}
\| KN_1(A_n)\|_{\mathcal{L}(L^2_{\ast}(\partial\Omega),L^2(\Omega,\mathbb{R}^N))}-1/n\leq\\ \|KN_1(A_n)[g_n]\|_{L^2(\Omega,\mathbb{R}^N)}
 \leq \| KN_1(A_n)\|_{\mathcal{L}(L^2_{\ast}(\partial\Omega),L^2(\Omega,\mathbb{R}^N))}.
 \end{multline*}

Let us consider a subsequence such that
\begin{multline*}
\lim_k
\| KN_1(A_{n_k})\|_{\mathcal{L}(L^2_{\ast}(\partial\Omega),L^2(\Omega,\mathbb{R}^N))}
=\\\limsup_n \| KN_1(A_n)\|_{\mathcal{L}(L^2_{\ast}(\partial\Omega),L^2(\Omega,\mathbb{R}^N))}
\end{multline*}
and such that $g_{n_k}$ converges, as $k\to\infty$, to $g$ strongly in $W^{-1/2,2}_{\ast}(\partial \Omega)$.
Then
\begin{multline*}
\|KN_1(A_{n_k})[g_{n_k}]\|_{L^2(\Omega,\mathbb{R}^N)}
 \leq \\
\|KN_1(A_{n_k})[g_{n_k}-g]\|_{L^2(\Omega,\mathbb{R}^N)}+\|KN_1(A_{n_k})[g]\|_{L^2(\Omega,\mathbb{R}^N)}
 \leq\\ C\|g_{n_k}-g\|_{W^{-1/2,2}_{\ast}(\partial \Omega)} +\|KN_1(A_{n_k})\|_{\mathcal{L}(L^2_{\ast}(\partial\Omega),L^2(\Omega,\mathbb{R}^N))}.
 \end{multline*}
Passing to the limit as $k\to \infty$ we obtain that
\begin{multline*}
\|KN_1(A)\|\leq \liminf_n \|KN_1(A_n)\|\leq 
 \limsup_n \|KN_1(A_n)\|=\\ \lim_k \|KN_1(A_{n_k})\|\leq
 \limsup_k (\|KN_1(A_{n_k})[g_{n_k}]\|_{L^2(\Omega,\mathbb{R}^N)}+1/n_k)\leq \|KN_1(A)\|
\end{multline*}
where $\|\cdot\|=\|\cdot\|_{\mathcal{L}(L^2_{\ast}(\partial\Omega),L^2(\Omega,\mathbb{R}^N))}$.
Thus also the proof of the second part is concluded.\cvd

\bigskip

Finally, we prove the statement of Example~\ref{nonexistence}.
First of all we need to go back to the uniqueness issue. Actually, for any $C^{1}$ diffeomorphism $\varphi=(\varphi_1,\varphi_2)$ of $\overline{\Omega}$ onto itself such that $\varphi|_{\partial\Omega}=Id$ and any 
conductivity tensor $A$ in $\Omega$ we may define the push-forward of the conductivity tensor $A$ by $\varphi$ as
$$\varphi_{\ast}(A)(y)=\frac{J(x)A(x)J(x)^T}{|\det J(x)|}\quad\text{for any }y\in\Omega$$
where $J(x)$ is the Jacobian matrix of $\varphi$ in $x$ and $x=\varphi^{-1}(y)$. It is a well-known fact that
$$D\text{-}N(\varphi_{\ast}(A))=D\text{-}N(A)\quad\text{and}\quad
N\text{-}D(\varphi_{\ast}(A))=N\text{-}D(A).$$

An interesting issue for the inverse conductivity problem in the anisotropic case is whether this invariance by diffeomorphisms that leave unchanged the boundary is the only obstruction to uniqueness, at least for symmetric conductivity tensors. This problem has been solved in dimension $2$, first in \cite{Syl} for smooth conductivity tensors, then in \cite{Ast e Pai e Las} for the general $L^{\infty}$ case where the following result has been proved.

\begin{teo} Let $\Omega\subset \mathbb{R}^2$ be a bounded, simply connected domain with Lipschitz boundary. Let us call $\mathcal{M}_{sym}$ the class of symmetric conductivity tensors in $\Omega$.
For any $A\in \mathcal{M}_{sym}$ let us define
\begin{multline*}
\Sigma(A)=\{A_1\in\mathcal{M}_{sym} :\ A_1=\varphi_{\ast}(A)\\\text{ where }\varphi:\Omega\to\Omega\text{ is a $W^{1,2}$ diffeomorphism and }\varphi|_{\partial\Omega}=Id\}.
\end{multline*}

Then the Dirichlet-to-Neumann map $D\text{-}N(A)$, or equivalently the Neumann-to-Dirichlet map 
 $N\text{-}D(A)$, uniquely determines the class $\Sigma(A)$.
 \end{teo}

Therefore if  $D\text{-}N(\sigma I_N)=
D\text{-}N(\hat{A})$, or equivalently $N\text{-}D(\sigma I_N)=
N\text{-}D(\hat{A})$,
then $\hat{A}=\varphi_{\ast}(\sigma I_N)$ for a suitable
$W^{1,2}$ diffeomorphism $\varphi$ such that $\varphi|_{\partial\Omega}=Id$.
The following example, that has been kindly communicated by Giovanni Alessandrini, shows that for some symmetric conductivity tensor $\hat{A}$ this is actually impossible, therefore existence of a minimizer in the class of scalar conductivities may fail.

\begin{exam}\label{Giovanni}
Let us consider $\Omega=B_1\subset \mathbb{R}^2$ and two positive constants $a$ and $b$, with $a\neq b$. Let $\hat{A}$ be as in \eqref{Giovannicond}.

We assume, by contradiction, that there exists a $W^{1,2}$ diffeomorphism $\varphi=(\varphi_1,\varphi_2): B_1\to B_1$, with
$\varphi|_{\partial B_1}=Id$, and a scalar conductivity $\sigma I_N$ such that
$$\hat{A}(y)=\frac{J(x)(\sigma(x) I_N)J(x)^{T}}{|\det J(x)|}\quad\text{for almost any }y\in B_1,$$
where $J(x)$ is the Jacobian matrix of
$\varphi$ in $x$ and $x=\varphi^{-1}(y)$.
We obtain that for any $x\in B_1$
$$\left[\begin{matrix}a& 0\\0 &b\end{matrix}\right]=\frac{\sqrt{ab}}{|\det J(x)|}\left[\begin{matrix}|\nabla \varphi_1(x)|^2& \nabla\varphi_1(x)\cdot\nabla\varphi_2(x)\\ \nabla\varphi_1(x)\cdot\nabla\varphi_2(x) &|\nabla \varphi_2(x)|^2\end{matrix}\right].$$
By equalling the determinants, we readily infer that $\sigma(x)=\sqrt{ab}$ for any $x\in B_1$. Furthermore, by the results in \cite{Syl} we deduce that $\varphi$ is actually much smoother, for instance it is a $C^1$ diffeomorphism.
We conclude that $\nabla \varphi_2(x)= \lambda(x) \left[\begin{smallmatrix}0 & -1\\1 & 0
\end{smallmatrix}\right]\nabla \varphi_1(x)$ with $\lambda(x)$ being either $\sqrt{b/a}$ or $-\sqrt{b/a}$. By the smoothness of $\varphi$ and the fact that $\varphi$ is a diffeomorphism, we deduce that $\lambda(x)$ is identically equal to $\sqrt{b/a}$ or to $-\sqrt{b/a}$ all over $B_1$.
We conclude that
$$\Delta \varphi_1=\Delta\varphi_2=0\quad\text{in }B_1.$$
Since $\varphi_1(x_1,x_2)=x_1$ and $\varphi_2(x_1,x_2)=x_2$ on $\partial B_1$, we obtain that
$\varphi$ is the identity and thus we have a contradiction.
\end{exam}

We are in the position of showing that Example~\ref{nonexistence} holds true.
In fact, by Proposition~\ref{approxbyscalar} and by the continuity properties of Theorem~\ref{Hcontinuity}, we have
$$\inf_{A\in \mathcal{M}_{scal}(\alpha_1,\beta_1)}
\|N\text{-}D(A)-\hat{\Lambda}_N\|=0
\quad\text{and}\quad \inf_{A\in \mathcal{M}_{scal}(\alpha_1,\beta_1)}
\|KN_1(A)\|=0.$$
Therefore a minimum does not exists unless there exists a scalar conductivity $\sigma I_N$ such that $\|N\text{-}D(\sigma I_N)-\hat{\Lambda}_N\|=0$ or $\|KN_1(\sigma I_N)\|=0$, that is $N\text{-}D(\sigma I_N)=
N\text{-}D(\hat{A})$, and this is impossible by Example~\ref{Giovanni}.

\appendix

\section{Appendix}
\label{Lpsection}

In this appendix, for completeness we treat the continuity with respect to $L^q$ norms.
At the end we provide as a consequence of these results a simple proof of Proposition~\ref{L1impliesH}.

We fix $\Omega\subset \mathbb{R}^N$, $N\geq 2$, a bounded domain with Lipschitz boundary.
The $L^{\infty}$ case is a classical and easy result that we state here.

\begin{teo}\label{inftycontinuity}
Let us consider two conductivity tensors $A_1$ and $A_2\in \tilde{\mathcal{M}}(\alpha,\tilde{\beta})$.

We have that $D\text{-}N$ and $N\text{-}D$ are Lipschitz continuous with respect to the $L^{\infty}(\Omega)$ norm on $\tilde{\mathcal{M}}(\alpha,\tilde{\beta})$ and the natural operator norms, that is
$$\|D\text{-}N(A_1)-D\text{-}N(A_2)\|_{\mathcal{L}(W^{1/2,2}(\partial \Omega),W^{-1/2,2}_{\ast}(\partial \Omega))}\leq C
\|A_1-A_2\|_{L^{\infty}(\Omega)}$$
and
$$\|N\text{-}D(A_1)-N\text{-}D(A_2)\|_{\mathcal{L}(W^{-1/2,2}_{\ast}(\partial \Omega),W^{1/2,2}_{\ast}(\partial \Omega))}\leq C
\|A_1-A_2\|_{L^{\infty}(\Omega)}$$
where $C$ depends on $N$, $\Omega$, $\alpha$ and $\tilde{\beta}$ only.

\end{teo}


Theorem~\ref{inftycontinuity} would suggest that the natural 
metric on
$\tilde{\mathcal{M}}(\alpha,\tilde{\beta})$ is the one induced by the $L^{\infty}$ norm. However, this is not always suited for applications, in particular in the case of 
discontinuous conductivity tensors. Therefore we consider  a weaker norm, namely the one induced by the $L^1$ norm (or, equivalently, given the uniform $L^{\infty}$ bound, by the $L^q$ norm for some $q$, $1\leq q<+\infty$). In order to obtain continuity, however, we need to make the following restriction.
We need to consider our Dirichlet-to-Neumann maps, or Neumann-to-Dirichlet maps respectively, as linear operators on a suitable Banach space $B_1$ which is continuously immersed  in $W^{1/2,2}(\partial\Omega)$, or on a suitable Banach space $\tilde{B}_1$ which is continuously immersed  in $W^{-1/2,2}_{\ast}(\partial\Omega)$ respectively,
with the corresponding norms. The main idea of these results goes back to \cite{Ron08}, where only the Neumann-to-Dirichlet case was treated.
The crucial remark, as already pointed out in \cite{Ron08}, is the use of a theorem by Meyers which we state here.

For any $p$, $1<p<+\infty$, let $W^{1-1/p,p}(\partial \Omega)$ be the space of traces of
$W^{1,p}(\Omega)$ functions on $\partial \Omega$.
We recall that for any $p$, $1\leq p\leq +\infty$, $p'$ denotes its conjugate exponent, that is $1/p+1/p'=1$.

\begin{prop}\label{regprop1}
Let $A\in \tilde{\mathcal{M}}(\alpha,\tilde{\beta})$.
There exists a constant $Q_1>2$, depending on $N$, $\Omega$, $\alpha$ and $\tilde{\beta}$ only, such that
for any $p$, $2<p<Q_1$, the following holds.

There exists a constant $C$, depending on $N$, $\Omega$, $\alpha$, $\tilde{\beta}$ and $p$ only,
such that for any 
$\varphi\in W^{1,p}(\Omega)$ and $F\in (W^{1,p'}(\Omega))^{\ast}$ we have
for $u$ solution to \eqref{Dirichletproblem}
\begin{equation}\label{MeyDir}
\|u\|_{W^{1,p}(\Omega)}\leq C(\|\varphi\|_{W^{1,p}(\Omega)}+\|F\|_{(W^{1,p'}\Omega))^{\ast}})
\end{equation}
and for $v$ solution to \eqref{Neumannproblem}, with $F$ such that $\langle F,1\rangle=0$,
\begin{equation}\label{MeyNeu}
\|v\|_{W^{1,p}(\Omega)}\leq C\|F\|_{(W^{1,p'}\Omega))^{\ast}}.
\end{equation}
\end{prop}


\proof{.} The Dirichlet case is treated in the original paper by Meyers, \cite{NMey}, whereas the extension to Neumann problems is contained in \cite{Gal e Mon}.\cvd

\bigskip

With the same notation as before, for any $A_1$, $A_2\in \tilde{\mathcal{M}}(\alpha,\tilde{\beta})$, for any $p$, $2<p<Q_1$, and any
$\varphi\in W^{1,p}(\Omega)$ and $F\in (W^{1,p'}(\Omega))^{\ast}$, we have
\begin{equation}\label{Direst}
\|u_1-u_2\|_{W^{1,2}(\Omega)}\leq C(\|\varphi\|_{W^{1,p}(\Omega)}+\|F\|_{(W^{1,p'}\Omega))^{\ast}})\|A_1-A_2\|_{L^q(\Omega)}
\end{equation}
and, provided $\langle F,1\rangle=0$, 
\begin{equation}\label{Neuest}
\|v_1-v_2\|_{W^{1,2}(\Omega)}\leq C\|F\|_{(W^{1,p'}\Omega))^{\ast}}\|A_1-A_2\|_{L^q(\Omega)}
\end{equation}
where $C$ depends on $N$, $\Omega$, $\alpha$, $\tilde{\beta}$ and $p$ only, and $q$,
$2<q<+\infty$, is such that
\begin{equation}\label{qdef}
\frac{1}{q}+\frac{1}{p}+\frac{1}{2}=1.
\end{equation}

Therefore we can easily deduce the following H\"older continuity result.

\begin{teo}\label{pcontinuity}
Let us consider two conductivity tensors $A_1$ and $A_2\in \tilde{\mathcal{M}}(\alpha,\tilde{\beta})$.


Fixed $p$, $2<p<Q_1$, $D\text{-}N$ and $N\text{-}D$ are H\"older continuous with respect to the $L^1(\Omega)$ norm on $\tilde{\mathcal{M}}(\alpha,\tilde{\beta})$ and the following operator norms.
We have
$$\|D\text{-}N(A_1)-D\text{-}N(A_2)\|_{\mathcal{L}(W^{1-1/p,p}(\partial \Omega),W^{-1/2,2}_{\ast}(\partial \Omega))}\leq C
\|A_1-A_2\|^{\delta}_{L^1(\Omega)}.$$

For any Banach spaces $\tilde{B}_1$ and $\tilde{B}_2$ such that $\tilde{B}_1$ is contained in the subspace of $g$ belonging to the
dual of $W^{1-1/p',p'}(\partial \Omega)$ such that $\langle g,1\rangle =1$, and $W^{1/2,2}_{\ast}(\partial \Omega)$ is contained in $\tilde{B}_2$, in both cases with continuous immersions, we have
$$\|N\text{-}D(A_1)-N\text{-}D(A_2)\|_{\mathcal{L}(B_1,B_2)}\leq C
\|A_1-A_2\|^{\delta}_{L^1(\Omega)}$$
where $C$ and $\delta$, $0<\delta<1$, depend on $N$, $\Omega$, $\alpha$, $\tilde{\beta}$ and $p$ only.
\end{teo}

A particularly interesting case for Neumann-to-Dirichlet maps is the following.
We can choose $\tilde{B}_1=\tilde{B}_2=L^2_{\ast}(\partial\Omega)$ since 
$L^2(\partial\Omega)$ is contained in the dual of $W^{1-1/p',p'}(\partial \Omega)$ for some $p$, $2<p<Q_1$, with $p$ close enough to $2$,
and $W^{1/2,2}_{\ast}(\partial\Omega)\subset L^2_{\ast}(\partial\Omega)$, with continuous immersions.

As a corollary to these results we can give a simple proof of Proposition~\ref{L1impliesH}.

\proof{ of Proposition~\ref{L1impliesH}}. First of all we assume that $\Omega$ has a Lipschitz boundary and that the convergence is in $L^1$ instead of $L^1_{loc}$.
Let us consider $F\in L^2(\Omega)$ and $\varphi=0$.
We denote $u_n=u(A_n,F,0)$ and $u=u(A,F,0)$ as before.

We notice that there exists a constant $p$, $2<p<Q_1$, such that $L^2(\Omega)\subset
(W^{1,p'}(\Omega))^{\ast}$ with continuous immersion. Therefore, by \eqref{Direst},
$u_n$ converges to $u$ strongly in $W^{1,2}(\Omega)$ as $n\to \infty$. About $A_n\nabla u_n$
we have that, up to a subsequence, it converges to $V$ weakly in $L^2(\Omega,\mathbb{R}^N)$. On the other hand
we notice that, by the uniform bound, $A_n$ converges to $A$ also in $L^2(\Omega)$, therefore
$A_n\nabla u_n$ converges to $A\nabla u$ strongly in $L^1(\Omega,\mathbb{R}^N)$, hence
$V=A\nabla u$ and the whole sequence $A_n\nabla u_n$ converges to $A\nabla u$ weakly in $L^2(\Omega,\mathbb{R}^N)$.

For the general case, we notice that, up to a subsequence, $A_n$ $H$-converges to $\tilde{A}$. By the locality of $H$-convergence and the previous part of the proof, we conclude that $\tilde{A}$ must coincide with $A$ almost everywhere. Therefore, the whole sequence $A_n$ $H$-converges to $A$.\cvd

\subsubsection*{Acknowledgements}
The author is partly supported by Universit\`a degli Studi di Trieste through Fondo per la Ricerca di Ateneo -- FRA 2012 and by GNAMPA, INdAM.
The author wishes to thank Adrian Nachman for stimulating discussions and Giovanni Alessandrini for  useful discussions and for kindly communicating Example~\ref{Giovanni}.

\end{document}